\documentclass[12pt,a4paper,twoside]{article}
\usepackage{a4wide}
\usepackage{amssymb}
\usepackage{amsmath}
\newtheorem{theorem}{Theorem}[section]

\newtheorem{corollary}[theorem]{Corollary}
\newtheorem{conjecture}[theorem]{Conjecture}

\newenvironment{pf}{\noindent{\bf Proof.}}{\hfill$\Box$}

\newenvironment{ac}{\noindent{\bf Acknowledgements.}}{}
\def\ii{{\mathbb I}}

\def\KK{{\mathbb K}}

\def\AA{{\mathcal A}}

\def\QQ{{\mathbb Q}}
\def\RR{{\mathbb R}}
\def\ZZ{{\mathbb Z}}

\begin{document}
\title{On the irreducibility of the determinant of the matrix of moving planes}
\author{Carlos D' Andrea \\ Departament of Mathematics \\ University of California \\94720 Berkeley, CA \\ USA \\ {\tt cdandrea@math.berkeley.edu}}
\date{}
\maketitle

\begin{abstract}
We show that, for generic bihomogeneous polynomials, the determinant of the matrix of moving planes is irreducible.
\end{abstract}

\section{Introduction}
Let $m$ and $n$ be positive integers, and
\begin{equation}\label{sup}
\AA:=\{(a_1,a_2)\in\ZZ^2: 0\leq a_1\leq m, 0\leq a_2\leq n\}.
\end{equation}
For $i=1,2,3,4,$ set $F_i:=\sum_{a\in\AA} c_{ia}x^a,$ where $c_{ia},x_1,x_2$ are indeterminates, and $x^a$ stands for
$x_1^{a_1}x_2^{a_2}$ if $a=(a_1,a_2).$
\par
Let $\KK:=\QQ(c_{ia})_{i=1,\dots,4,\,a\in\AA},$ i.e. the field of quotients of the indeterminates $c_{ia},$
and denote with $S_{i,j}$ the $\KK$-vector space generated by those
monomials whose degree in $x_1$ (resp.\ $x_2$) is less than or equal to
$i$ (resp.\ $j$).
Consider the following $\KK$-linear map
\begin{equation}
\label{map}
\begin{array}{ccc}
{S_{m-1,n-1}}^4&\stackrel{\phi_{m,n}}{\to}&S_{S_{2m-1,2n-1}}\\
(A_1,A_2,A_3,A_4)&\mapsto&\sum_{i=1}^4 A_iF_i,
\end{array}
\end{equation}
and denote with $F_{m,n}$ its matrix in the monomial bases of ${S_{m-1,n-1}}^4$ and $S_{2m-1,2n-1}.$ It is a square matrix of size
$4mn\times 4mn$ whose entries are actually in $\ZZ[c_{ia}].$ Note that $\det(F_{m,n})$ is well-defined up to a sign.
\par The main result of this note is the following
\begin{theorem}\label{irred}
$\det(F_{m,n})$ is an irreducible element of $\ZZ[c_{ia}].$
\end{theorem}

In order to set this result in a proper context, we give here a brief introduction to the method of moving quadrics for the implicitization of rational
surfaces \cite{SC,CGZ,BCD,DK}. For any specialization of the $c_{ia}$ in any  field ${\bf K},$ we can think of the quadruple $(F_1,F_2,F_3,F_4)$ as
the parametrization of the following rational surface in ${\bf K}^3:$
\begin{equation}\label{imp}
\begin{array}{ccc}
Y_1(x)=\frac{F_1(x)}{F_4(x)}&Y_2=\frac{F_2(x)}{F_4(x)}&Y_3(x)=\frac{F_3(x)}{F_4(x)}.
\end{array}
\end{equation}
The implicitization problem consists in the computation of the
irreducible polynomial $P(Y_1,Y_2,Y_3)$ whose vanishing defines the closure of
this surface.
\par In \cite{SC}, Sederberg and Chen introduced a new technique for finding
the implicit equation called the \textit{method of moving quadrics}.
This method is based in the construction of a matrix whose entries are the coefficients in the monomial bases of certain syzygies of $I^2,$
where $I$ is the ideal generated by $F_1,F_2,F_3,F_4.$
\par In \cite[Theorem $4.2$]{CGZ}, it is shown that the method works in this case (tensor product surfaces) provided that
\begin{itemize}
\item ${\rm Res}_{m,n}(F_1,F_2,F_3)\neq0,$ where ${\rm Res}_{m,n}$ stands for the sparse resultant as defined in \cite{CLO}, associated with the
support $\AA$ defined in (\ref{sup}).
\item $\det(F_{m,n})\neq0.$
\end{itemize}
Let $Q_{m,n}$ be the matrix in the monomial bases of the following map
\begin{equation}\label{quad}
{S_{m-1,n-1}}^9\to S_{3m-1,3n-1}
\end{equation}
which sends $(A_1,A_2,A_3,A_4,A_5,A_6,A_7,A_8,A_9)$ to
{\small
$$A_1F_1^2+A_2F_2^2+A_3F_3^2+A_4F_1F_2+A_5F_1F_3+A_6F_1F_4+A_7F_2F_3+A_8F_2F_4+A_9F_3F_4.$$
}
In \cite{CGZ} it is also shown that if the method of the moving quadrics works, then $\det(Q_{m,n})\neq0,$ and it was conjectured that
\begin{equation}\label{eqq}
\det(Q_{m,n})=\pm{\det(F_{m,n})}^3{\rm Res}_{m,n}(F_1,F_2,F_3)
\end{equation}
(see \cite[Section $6$]{CGZ}). This conjecture was based on a similar relationship holding in the case of rational
curves (\cite{SGD,ZCG}), and it was hinted that the proof of this equality would require to show that $\det(F_{m,n})$ is
an irreducible polynomial as in the curve case.
\par In \cite{dan}, we have proven that (\ref{eqq}) holds by using homological algebra methods, like the Cayley method for computing
the determinant of a generically exact Koszul complex. The irreducibility of the determinant of $F_{m,n}$ was not needed.
\par Theorem \ref{irred} implies that the right hand side of (\ref{eqq})
is the factorization of $\det(F_{m,n})$ into irreducible components in $\ZZ[c_{ia}].$ The complete factorization of other
determinants arising from the general version of (\ref{eqq}) given in \cite{dan} also holds from Theorem \ref{irred}.
\par
The fact that these determinants are irreducible may provide a geometric meaning to them as in the case of resultants, whose irreducibility is
explained because they are projection operators from some irreducible variety.
Recently, the irreducibility of some subresultants have been proven in \cite{BD} by
using geometric methods, but it should be noted that although the determinant of $F_{m,n}$ may be regarded as an \textit{inertia form} of the
ideal $I$ in the sense of Jouanolou (\cite{jou1}), it cannot be considered as a \textit{multivariate subresultant} as defined by Chardin in \cite{cha} due to the
fact that in this case the number of polynomials is greater than the number of (homogeneous) variables and so the theory of multivariate subresultants cannot
be applied.
\par
The paper is organized as follows. In Section \ref{dos}, we give the proof of Theorem \ref{irred}.
This proof is rather elementary and follows the same ideas that led to the proofs given by Macaulay in \cite{Mac} for the irreducibility of
the Sylvester resultant and by Sederberg, Goldman and Du in \cite{SGD} for the irreducibility of the determinant of the matrix of moving lines, the equivalent of $F_{m,n}$ in
the univariate case.
Then, we show as a corollary of Theorem \ref{irred} that (\ref{eqq}) holds without applying the results of \cite{dan}.
The paper concludes with the statement of a general conjecture in Section \ref{quatro}.

\smallskip
\begin{ac}
I am grateful to Ron Goldman for helpful discussions and to David Cox for comments on preliminary versions of this draft.
\end{ac}

\bigskip
\section{Proof of Theorem \ref{irred}}\label{dos}
We shall prove Theorem \ref{irred} by double induction on $m$ and $n.$
\subsection{Case $n=1$}
First, we shall see that Theorem \ref{irred} holds for $n=1$ and any $m.$ In order to do this, we will make induction on $m.$
\par For $m=1,$ we have that
$$F_{1,1}=\left(\begin{array}{cccc}
c_{1,(0,0)}& c_{1,(1,0)}&c_{1,(0,1)}&c_{1,(1,1)}\\
c_{2,(0,0)}& c_{2,(1,0)}&c_{2,(0,1)}&c_{2,(1,1)}\\
c_{3,(0,0)}& c_{3,(1,0)}&c_{3,(0,1)}&c_{3,(1,1)}\\
c_{4,(0,0)}& c_{4,(1,0)}&c_{4,(0,1)}&c_{4,(1,1)}
\end{array}\right).$$
As $F_{1,1}$ is a generic matrix of $4\times4,$ its determinant is clearly irreducible.
\par Now take $m>1$ and suppose that $\det(F_{m,1})=P(c)*Q(c),$ where $P,Q\in\QQ[c_{ia}].$
As $\det(F_{m,1})$ is a homogeneous polynomial, then $P$ and $Q$ must also be homogeneous.
\par
Let
$$\tilde{c}_{ia}:=\left\{\begin{array}{cc}
c_{ia}&\mbox{if}\,a_1=0\,\mbox{or}\,a_1=m\\
0&\mbox{otherwise,}
\end{array}
\right.$$
and denote with $\Delta$ the determinant of
$$\left(\begin{array}{cccc}
c_{1,(0,0)}& c_{1,(m,0)}&c_{1,(0,1)}&c_{1,(m,1)}\\
c_{2,(0,0)}& c_{2,(m,0)}&c_{2,(0,1)}&c_{2,(m,1)}\\
c_{3,(0,0)}& c_{3,(m,0)}&c_{3,(0,1)}&c_{3,(m,1)}\\
c_{4,(0,0)}& c_{4,(m,0)}&c_{4,(0,1)}&c_{4,(m,1)}
\end{array}\right).$$
It is straigthforward to check that $\det(F_{m,1}(\tilde{c}))=\pm\Delta^m.$ This shows that $\det(F_{m,1})$ is not identically zero in
$\KK,$ and moreover it has content $1$ in $\ZZ[c_{ia}].$ So, in order to prove the irreducibility, it is enough if we show that
$\det(F_{m,1})$ is irreducible in $\QQ[c_{ia}].$
\par
As $\Delta$ is irreducible, there exists $k,\,0\leq k\leq m$ such that
$\deg(P)=4k$ and $\deg(Q)=4(m-k).$
\par In addition, if we set $\omega(c_{ia}):=a_2,$ then it is straightforward to check that
$\det(F_{m,1})$ is homogeneous with respect to the weight given by $\omega,$ and its ``$\omega$-degree'' is $2m^2.$
This notion of homogeneity with respect to a weight has already been used by Macaulay in \cite{Mac} with the name
``isobarism''.
\par
So, $P$ and $Q$ must be also $\omega$-homogeneous, and as $\Delta$ has $\omega$-degree equal to $2m,$ then $P$ must have $\omega$-degree equal to
$2mk,$ and $Q$ must have $2(m-k)m.$
\par
Consider now another specialization
$$\overline{c}_{ia}:=\left\{\begin{array}{cc}
t*c_{ia}&\mbox{if}\,a_2=0,\\
c_{ia}&\mbox{otherwise.}
\end{array}
\right.$$
As in \cite{Mac}, we shall index all matrices in such a way that the rows correspond to the elements of the monomial basis of the
domain of the transformation.
If we sort the columns of $F_{m,n}$ as follows: $\{1,x_1,x_2,x_1*x_2,\dots\},$ then after the last specialization,
$F_{m,n}$ will be
\begin{equation}
\label{tul}
\left(
\begin{array}{ccccccc}
tc_{1,(0,0)}& tc_{1,(1,0)}&c_{1,(0,1)}&c_{1,(1,1)}&\dots&\dots&\dots\\
tc_{2,(0,0)}& tc_{2,(1,0)}&c_{2,(0,1)}&c_{2,(1,1)}&\dots&\dots&\dots\\
tc_{3,(0,0)}& tc_{3,(1,0)}&c_{3,(0,1)}&c_{3,(1,1)}&\dots&\dots&\dots\\
tc_{4,(0,0)}& tc_{4,(1,0)}&c_{4,(0,1)}&c_{4,(1,1)}&\dots&\dots&\dots\\
0&0&tc_{1,(0,0)}& tc_{1,(1,0)}&\star&\star&\star\\
0&0&tc_{2,(0,0)}& tc_{2,(1,0)}&\star&\star&\star\\
0&0&tc_{3,(0,0)}& tc_{3,(1,0)}&\star&\star&\star\\
0&0&tc_{4,(0,0)}& tc_{4,(1,0)}&\star&\star&\star\\
0&0&0&0&\star&\star&\star\\
\vdots&\vdots&\vdots&\vdots&\star&\star&\star\\
0&0&0&0&\star&\star&\star
\end{array}\right),
\end{equation}
where the $``\star''$ area is a square matrix of size
$4(m-1)\times4(m-1)$ corresponding to a matrix of the form $F_{m-1,1}.$
We call this matrix $F'_{m-1,1}$ because the coefficients are generic but indexed in
a different way. Observe also that $\det(F'_{m-1,n})$ is $\omega$-homogeneous, but its $\omega$-degree is not $2(m-1)^2.$
\par
Let $\overline{\Delta}$ be the determinant of
$$\left(\begin{array}{cccc}
c_{1,(0,0)}& c_{1,(1,0)}&c_{1,(0,1)}&c_{1,(1,1)}\\
c_{2,(0,0)}& c_{2,(1,0)}&c_{2,(0,1)}&c_{2,(1,1)}\\
c_{3,(0,0)}& c_{3,(1,0)}&c_{3,(0,1)}&c_{3,(1,1)}\\
c_{4,(0,0)}& c_{4,(1,0)}&c_{4,(0,1)}&c_{4,(1,1)}
\end{array}\right).$$
So, we have that
$$
\det(F_{m,1}(\overline{c}))=t^2\overline{\Delta}\det(F'_{m-1,1})+
\,\mbox{higher order terms in},\, t.$$
From now on, the ``degree'' $\deg$  of a polynomial in the variables $c_{ia}$ will be a pair
of the form $($homogeneous degree, $\omega$-degree $).$ So, we have the following:
\begin{itemize}
\item $\deg(\overline{\Delta})=(4,2).$
\item $\deg(\det(F'_{m-1,1}))=(4(m-1),2m^2-2).$
\end{itemize}
As $\det(F_{m,1}(\overline{c}))=P(\overline{c})*Q(\overline{c})$ in
$\QQ[c_{ia}][t],$ and due to the fact that  $\det(F'_{m-1,1})$ is irreducible by the inductive hypothesis,
this determinant must be a factor of the lowest term in $t$ of $P(\overline{c})$ or
$Q(\overline{c}).$ So, we have the following two scenarios:
\begin{enumerate}
\item If it is a factor of $P,$ then $(4(m-1),2m^2-2)\leq(4k,2mk).$ But this implies
that $m=k$ and so $Q$ is a constant.
\item If it is a factor of $Q,$ then $(4(m-1),2m^2-2)\leq(4(m-k),2m(m-k)),$ and this
inequality implies that $k=0,$ and so $P$ is a constant.
\end{enumerate}
So, $\det(F_{m,1})$ is irreducible and we are done.
\par
\bigskip
\subsection{The General Case}
Now we are going to give the proof for the general case.
We fix $n$ and will make induction on $m=1,2,\dots,n.$
The initial case has already been proved in the previous section.
\par Take now $1<m\leq n,$ and suppose as before that $\det(F_{m,n})=P(c)*Q(c).$
Setting $\omega(c_{ia}):=a_1,$ it is straightforward to check that $\det(F_{m,n})$ is
$\omega$-homogeneous, and its degree is $(4mn,2m^2n)$ (here we keep the notation of the previous section, i.e. ``degree'' is (homogeneous degree,
$\omega$-degree) ).
\par
Consider the following specialization:
Let
$$\tilde{c}_{ia}:=\left\{\begin{array}{cc}
c_{ia}&\mbox{if}\,a\in\{(0,0),(0,n),(m,0),(m,n)\}\\
0&\mbox{otherwise,}
\end{array}
\right.$$
and denote with $\Delta$ the determinant of
$$\left(\begin{array}{cccc}
c_{1,(0,0)}& c_{1,(m,0)}&c_{1,(0,n)}&c_{1,(m,n)}\\
c_{2,(0,0)}& c_{2,(m,0)}&c_{2,(0,n)}&c_{2,(m,n)}\\
c_{3,(0,0)}& c_{3,(m,0)}&c_{3,(0,n)}&c_{3,(m,n)}\\
c_{4,(0,0)}& c_{4,(m,0)}&c_{4,(0,n)}&c_{4,(m,n)}
\end{array}\right).$$
As $\Delta$ is irreducible and homogeneous of degree $(4,2m),$
and due to the fact that
$\det(F_{m,n}(\tilde{c}))=\Delta^{mn}=P(\tilde{c})*Q(\tilde{c}),$
it turns out that there must exist $k,\,0\leq k\leq mn-k$ such that
$P$ (resp. $Q$) is homogeneous of degree $(4k,2mk)$ (resp.
$(4(mn-k),2m(mn-k))$ ). This also shows that the determinant of $F_{m,n}$ is not identically zero in $\KK$ and has content one in
$\ZZ[c_{ia}],$ so in order to prove the claim it is enough to show that the irreducibility holds in $\QQ[c_{ia}].$
\par
Consider now the other specialization
$$\overline{c}_{ia}:=\left\{\begin{array}{cc}
t*c_{ia}&\mbox{if}\,a_1=0,\\
c_{ia}&\mbox{otherwise.}
\end{array}
\right.$$
If we index the columns of $F_{m,n}$ in such a way that the first $4m$ columns
are indexed by $x_1^{a_1}x_2^{a_2},$ with $a_1\leq 1,$ then it turns out that the
specialized matrix $F_{m,n}(\overline{c})$ will have a structure similar to
(\ref{tul}). To be more precise, it will have the following block structure:
$$\left(
\begin{array}{ccc}
A&B&\dots\\
0&C&F'_{m-1,n}
\end{array}\right),$$
Where each coefficient in $A$ and $C$ is a multiple of $t$
and the block $(A,B)$ if we set $t=1$ is $F_{m,1}.$
\par So, we have that
$$\det(F_{m,n}(\overline{c}))=t^{2m}\det(F_{m,1})\det(F'_{m-1,n})+
\mbox{higher terms in}\, t.$$
An explicit computation shows that
\begin{itemize}
\item $\deg(\det(F_{m,1}))=(4m,2m).$
\item $\deg(\det(F'_{m-1,n}))=(4mn-4m,2m^2n-2m).$
\end{itemize}
As before, if we regard $P(\overline{c})$ and $Q(\overline{c})$ as polynomials in $t,$
and due to the fact that both $\det(F'_{m-1,n})$ and $\det(F_{m,1})$ are irreducible
polynomials (here we use the
inductive hypothesis), it turns out that one of them must be contained in
$P(\overline{c})$ and the other must be contained in $Q(\overline{c})$ (otherwise,
either $P$ or $Q$ are constants are we are done).
\par Again, we have two different scenarios:
\begin{enumerate}
\item $$\left\{\begin{array}{l}
(4m,2m)=(4k,2mk),\\
(4mn-4m,2m^2n-2m)=(4(mn-k),2m(mn-k)).
\end{array}\right.
$$ This implies that $m=k=1$ which is impossible, since we are assuming $m>1.$
\item $$\left\{\begin{array}{l}
(4m,2m)=(4(mn-k),2m(mn-k))\\
(4mn-4m,2m^2n-2m)=(4k,2mk).
\end{array}\right.$$
This implies that $m=mn-k=1,$ again a contradiction.
\end{enumerate}
So, there cannot be such a factorization and hence $\det(F_{m,n})$ is irreducible.

\bigskip
\subsection{A complete factorization of the determinant of the moving quadric matrix}
A straightforward consequence of Theorem \ref{irred} is the following:
\begin{corollary}
$$\det(Q_{m,n})=\pm{\det(F_{m,n})}^3{\rm Res}_{m,n}(F_1,F_2,F_3).$$
\end{corollary}
\begin{pf}
In \cite[Theorem $4.1$]{CGZ}, it is shown that if we specialize the coefficients $c_{ia}$ in such a way that
${\rm Res}_{m,n}(F_1,F_2,F_3)\neq0$ and $\det(F_{m,n})\neq0,$ then $\det(Q_{m,n})\neq0.$
Then, Hilbert's Nullstellensatz and the irreducibility of the resultant and the determinant of $F_{m,n}$ implies that
there exist a constant ${\bf c}\in\QQ,$ and positive integers $a$ and $b$ such that
\begin{equation}\label{fi}
\det(Q_{m,n})={\bf c}\,{\det(F_{m,n})}^a{{\rm Res}_{m,n}(F_1,F_2,F_3)}^b.
\end{equation}
Computing the degree in the coefficients of $F_4$ in both sides of this equality, we get $3mn=amn,$ and from here we have that $a=3.$
In order to get $b=1,$ we compute the degree in the coefficients of $F_1$ in both sides of (\ref{fi}).
\par
Last, by expanding both sides of (\ref{fi}) as a polynomial in one of the extremal coefficients of $F_4$ (like $c_{4,(m,n)}$) and comparing the
leading terms, we will get that ${\bf c}=\pm1.$
\end{pf}

\section{The general conjecture}\label{quatro}
The method of moving quadrics was initially designed for the implicitization of surfaces given by triangular or tensor product parametrizations (see
\cite{CGZ}). Recently, it has been extended in \cite{DK} to parametrizations with fixed support. In order to state our general conjecture,
we will briefly review the results of \cite{DK}.
\par Let $\AA$ be any subset of $\ZZ^2,$ and set $Q$ the convex hull of $\AA.$
Denote with $E(t)$ the Ehrhart polynomial of $Q$ as defined in \cite{sta}. Write
$E(t)=At^2+\frac{B}{2}t+1,$
where $A=\frac{Area(Q)}{2}$ and $B$ equals the number of boundary points.
Let $E_I$ be a connected set of edges of $Q,$ we will denote with  $B_I$ the sum of the edge lenghts of $E_I.$
It turns out that the Ehrhart polynomial of $E_I$ equals $B_I t +1.$
\par Take $E_I$ in such a way that $B\geq 2B_I.$ Set as before $F_i:=\sum_{a\in\AA}c_{ia}x^a$ and $\KK:=\QQ(c_{ia}).$
For a set $P\subset\RR^2$ we will denote with $S_P$ the $\KK$ vector space generated by $x^a,$ with $a\in P\cap\ZZ^2.$
\par
Consider the following $\KK$-linear map
\begin{equation}\label{MP}
\begin{array}{cccc}
\phi_\AA:&{S_{Q\setminus E_I}}^4&\to&S_{2Q\setminus 2E_I} \\
&(A_1,A_2,A_3,A_4)&\mapsto&\sum_{i=1}^4 A_i F_i,
\end{array}
\end{equation}
and let $F_\AA$ be the matrix of this map in the monomial bases.

\medskip
Let $\ii\subset\left(Q\setminus E_I\right)\cap\ZZ^2$ such that $\# \ii=B-2B_I.$
Denote with $F_{\AA,\ii}$ the square submatrix of $F_\AA$ where we have deleted all the rows corresponding to $x^aF_4,\,a\in \ii.$

\smallskip
Consider now the following map
\begin{equation}\label{mq}
\begin{array}{cccl}
\psi_\AA:&{S_{Q\setminus E_I}}^{10}&\to&S_{3Q\setminus 3E_I} \\
&(A_{i,j,k,l})_{i+j+k+l=2}&\mapsto&\sum_{i+j+k+l=2}A_{i,j,k,l}F_1^iF_2^jF_3^kF_4^l,
\end{array}
\end{equation}
and let $Q_\AA$ be the matrix of $\psi_\AA$ in the monomial bases.

Following \cite{CGZ}, we remove all the rows indexed by $x^aF_4^2,\, a\in Q\setminus E_I,$ as well as the $3(B-B_I)$ rows
indexed by $x^a(F_1F_4,F_2F_4,F_3F_4),\,a\in \ii,$
and define $Q_{\AA,\ii}$ to be the coefficient matrix of the remaining polynomials.
One can check that $Q_{\AA,\ii}$ is a square matrix of order $9A+\frac32B-3B_I.$
\bigskip

The following is one of the main results of \cite{DK} and it is a generalization of the results given in \cite{CGZ} for the triangular and the
tensor product case:
\begin{theorem}\label{mt}
$$\det(Q_{\AA,\ii})={\det(F_{\AA,\ii})}^3 {\rm Res}_\AA(F_1,F_2,F_3).$$
\end{theorem}

The following conjecture  has been motivated by experimental
evidence, and should be regarded as the general version of Theorem
\ref{irred}:
\begin{conjecture}\label{kk}
Theorem \ref{mt} gives a complete factorization of $\det(Q_{\AA,\ii}),$ i.e. $\det(F_{\AA,\ii})$ is either zero or an irreducible element of
$\QQ[c_{ia}].$
\end{conjecture}

\end{document}